\documentclass[12pt]{article}

\usepackage{amsmath,amsthm,amsfonts,amssymb,latexsym,amscd}

\input{prepictex.tex}
\input{pictex.tex}
\input{postpictex.tex}

\begin{document}

\newtheorem*{theo}{Theorem}
\newtheorem*{pro}{Proposition}
\newtheorem*{cor}{Corollary}
\newtheorem*{lem}{Lemma}
\newtheorem{theorem}{Theorem}[section]
\newtheorem{corollary}[theorem]{Corollary}
\newtheorem{lemma}[theorem]{Lemma}
\newtheorem{proposition}[theorem]{Proposition}
\newtheorem{conjecture}[theorem]{Conjecture}
\newtheorem{definition}[theorem]{Definition}
\newtheorem{problem}[theorem]{Problem}
\newtheorem{remark}[theorem]{Remark}
\newtheorem{example}[theorem]{Example}
\newcommand{\Naturali}{{\mathbb{N}}}
\newcommand{\Reali}{{\mathbb{R}}}
\newcommand{\Complessi}{{\mathbb{C}}}
\newcommand{\Toro}{{\mathbb{T}}}
\newcommand{\Relativi}{{\mathbb{Z}}}
\newcommand{\HH}{\mathfrak H}
\newcommand{\KK}{\mathfrak K}
\newcommand{\LL}{\mathfrak L}
\newcommand{\as}{\ast_{\sigma}}
\newcommand{\tn}{\vert\hspace{-.3mm}\vert\hspace{-.3mm}\vert}
\def\mA{{\mathfrak A}}
\def\A{{\mathcal A}}
\def\mB{{\mathfrak B}}
\def\B{{\mathcal B}}
\def\C{{\mathcal C}}
\def\D{{\mathcal D}}
\def\F{{\mathcal F}}
\def\H{{\mathcal H}}
\def\J{{\mathcal J}}
\def\K{{\mathcal K}}
\def\L{{\cal L}}
\def\N{{\cal N}}
\def\M{{\cal M}}
\def\O{{\mathcal O}}
\def\P{{\cal P}}
\def\S{{\cal S}}
\def\T{{\cal T}}
\def\U{{\cal U}}
\def\W{{\cal W}}
\def\b{\lambda_B(P}
\def\j{\lambda_J(P}

\title{Computing the Jones index of quadratic permutation endomorphisms of $\O_2$}

\author{Roberto Conti, Wojciech Szyma{\'n}ski}

\date{}
\maketitle

\renewcommand{\sectionmark}[1]{}

\noindent
{\small \date{3 December, 2008}}

\vspace{7mm}
\begin{abstract}
We compute the index of the inclusions of type $III_{1/2}$ factors arising from
endomorphisms of the Cuntz algebra ${\mathcal O}_2$
associated to the rank-two permutation matrices.
\end{abstract}

\vfill
\noindent {\bf MSC 2000}: 46L37, 46L05

\vspace{3mm}
\noindent {\bf Keywords}: Cuntz algebra, endomorphism, automorphism,
 permutation, index, entropy.

\newpage

\section{Introduction}
Subfactor theory started with the breakthrough work by V. Jones \cite{Jo}.
In particular, given a subfactor $N \subset M$, i.e. an inclusion of von Neumann algebras
with trivial center, it is possible to associate a numerical invariant $[M:N]$ called
the index of $N$ in $M$. Roughly this invariant resembles the index of a subgroup, and it reduces to that
in very special situations, however in general it takes values
in the Jones set $J := \{4 \cos^2(\pi/n) \ | \ n=3,4,5,\ldots\} \cup [4,\infty]$.
The theory was developed first for subfactors of type $II_1$,
i.e. those admitting a faithful tracial state, but it has been generalized since to general
inclusions of factors by Kosaki and Longo \cite{Ko,Lo1}. In this more general framework the index
$Ind_E(N \subset M)$ depends on the choice of a faithful conditional expectation $E: M \to N$.
In turn, one can always find a unique conditional expectation that minimizes the index.
Subfactor theory has found important applications in Quantum Field Theory and Conformal Field Theory
(see e.g. \cite{Lo2,KaLo} for a small sample)
and in the study of the tensor categories arising from (compact) groups and quantum groups.
The square root of the minimal index is usually called statistical dimension, as it gives back
the DHR dimension of sectors in QFT/CFT and the quantum dimension in the context
of the representation theory of quantum groups.
When $M$ is an infinite factor the dimension provides a semi-ring homomorphism ${\rm Sect}_0(M) \to \sqrt{J}$,\footnote{That is,
additive for the natural direct sum of sectors and multiplicative for their composition.}
where ${\rm Sect}_0(M) := {\rm End}_0(M)/{\rm Inn}(M)$ and ${\rm End}_0(M)$ denotes the family of unital normal
$*$-endomorphisms of $M$ with finite index, i.e. those for which ${\rm Ind}(\rho(M) \subset M) < \infty$.
It has been also speculated that the Jones index could replace the Fredholm index in a truly
quantum form of the celebrated Atiyah-Singer theory \cite{l2}.
Anyhow, computing indices in specific situations is often challenging and sometimes worthwhile.
For a non-technical but broad overview of the subject including lot of important connections with other areas we refer the reader to \cite{j0}.

\medskip

For a number of years already, endomorphisms of the Cuntz algebras $\O_n$ have been subject
of intensive investigations. In addition to the intrinsic value of such study,
the reason for this wide interest in endomorphisms of $\O_n$ might be twofold.
Firstly, they provide a nice framework in which explicit computations
and analysis of non-trivial examples is feasible.
Secondly, they link in many an interesting way with several aspects of theory of operator algebras
such as subfactors and index theory, entropy, dynamical systems, wavelets and quantum field theory.

A very special class of such endomorphisms are the so-called permutation endomorphisms,
i.e. those arising from permutation unitaries in the core UHF-subalgebra.
In particular,
Kawamura \cite{Kaw1,Kaw2} analyzed systematically some properties
like properness, irreducibility and
branching laws of
rank-two (or ``quadratic'') permutation endomorphisms of $\O_2$.
More recently Skalski and Zacharias computed Voiculescu's topological entropy \cite{Vo}
of the same endomorphisms \cite{SZ}, thus extending Choda's computation of the entropy
of the canonical shift of $\O_n$ \cite{Ch}.

The main purpose of this
note is to add another bit of information,
namely the computation of the index in the sense of Jones-Kosaki-Longo \cite{Jo,Ko,Lo1,Lo2}
of the normal extension of these endomorphisms.
More precisely, we will mostly focus on the $III_{1/2}$ factors obtained as
the weak closure of $\O_2$ in the GNS representation of the
canonical KMS state with respect to the (rescaled) gauge action.
After Jones posed the problem in \cite{j},
general methods for computing the index of an inclusion of factors associated to a localized endomorphism
of a Cuntz algebra $\O_n$ were first foreseen by Longo \cite{l} using the theory of sectors.
Later such methods were developed in full strength in \cite{CP} for the type $III$
situation and in the unpublished PhD dissertation
of Akemann for the type $II_1$ case \cite{Ak} (cf. \cite[Example 5.1.6]{JS}),
see also \cite{CF,Go}.

\section{Main results}
Let ${\mathcal O}_2=C^*(s_1,s_2)$ be the Cuntz algebra \cite{Cun1}, that is the $C^*-$algebra generated by two isometries $s_1$, $s_2$ such that $s_1 s_1^* + s_2 s_2^* = 1$. Recall \cite{Cun2}
    that for any unitary $u\in{\mathcal O}_2$ there exists a unital $*$-endomorphism $\lambda_u$
    of ${\mathcal O}_2$ such that $\lambda_u(s_i) = u s_i$, $i=1,2$.
    Also, we denote by $\omega$ and $\pi_\omega$ the faithful canonical KMS state $\tau \circ E$ on $\O_2$
    and its GNS representation, respectively,
    where $E$ is the conditional expectation of $\O_2$ onto its core UHF subalgebra
    $\F_2 \simeq \otimes_{i=1}^\infty M_2$
    obtained by averaging over the gauge action of $\mathbb T$ and $\tau$ is the unique trace on $\F_2$.

In the case of rank-two endomorphisms $\lambda_u$ of $\O_2$
(that is,  when $u \in M_2 \otimes M_2 \subset \F_2$)
the minimal index ${\rm Ind}(\lambda_u)$ \cite{Hi}
of the inclusions of A.F.D. $III_{1/2}$ factors $\pi_\omega(\lambda_u(\O_2))'' \subset \pi_\omega(\O_2)''$
satisfies \cite{l,CP}
$${\rm Ind}(\lambda_u) \leq \big[\pi_\omega(\F_2)'': \pi_\omega(\lambda_u(\F_2))''\big] \leq 2^2 = 4$$
and thus all inclusions are extremal and ${\rm Ind}(\lambda_u)$ coincides with the Jones index
of the inclusion of A.F.D. $II_1$ factors $\pi_\omega(\lambda_u(\F_2))'' \subset \pi_\omega(\F_2)''$.
Moreover, according to the general analysis in \cite{CP,Ak},
in the above situation ${\rm Ind}(\lambda_u)$
is necessarily an integer and more precisely takes values in the set $\{1,2,4\}$.
The main result of this section is the computation of the precise values of the indices ${\rm Ind}(\rho_\sigma)$
for all the endomorphisms $\rho_\sigma := \lambda_{u_\sigma}$
of $\O_2$ induced by permutation matrices $u_\sigma \in M_2 \otimes M_2 \subset \F_2$.\footnote{
A few words on the notation: in the situation under consideration
these matrices $u_\sigma$ are naturally parametrized by permutations of the set $\{1,2\}^2$; after identifying $\{11,12,21,22\}$ with $\{1,2,3,4\}$, they are conveniently indexed by elements in the usual symmetric group on four elements.
Thus for instance $\rho_{12}$ is the endomorphism of $\O_2$ induced by the permutation matrix associated to the cycle $\sigma_{12} \equiv (12) \in {\mathfrak S}_4$.}
These values are plotted in the last column of Table 1 (below),
modelled on those in \cite{Kaw1,Kaw2,SZ}, where
for multiindices $\alpha$ and $\beta$
we denote $s_\alpha s^*_\beta$ by $s_{\alpha,\beta}$.\footnote{If $\alpha = (\alpha_1,\ldots,\alpha_k) \in \{1,2\}^k$
then $s_\alpha := s_{\alpha_1}\ldots s_{\alpha_k}$.}
Notice that by the results in \cite{Kaw1,Kaw2} among these 24 endomorphisms there are only 16
distinct inner equivalence classes.
Earlier computations included the well-known canonical endomorphism of $\O_2$, here denoted $\rho_{23}$,
and also $\rho_{12}$ \cite[Theorem 1.3 (ii)]{Kaw0}.
\medskip

\noindent{\bf Automorphisms of $\O_2$.}
There are four automorphisms of $\O_2$ \cite{Kaw1},
namely $\rho_{\rm id}$, $\rho_{(12)(34)}$, $\rho_{(13)(24)}$ and $\rho_{(14)(23)}$,
which clearly give index 1.

\medskip
\noindent{\bf Reducible endomorphisms of $\O_2$.}
There are ten reducible endomorphisms of $\O_2$ \cite{Kaw1}.
As
$\pi_\omega(\lambda_u(\O_2)' \cap \O_2) \subset \pi_\omega(\lambda_u(\O_2))' \cap \pi_\omega(\O_2)''$,
the corresponding inclusions of factors are reducible too and the index is necessarily $\geq (1+1)^2 = 4$.
Thus for all these endomorphisms the index is 4.

\medskip
In the sequel we follow closely the discussion in \cite[Section 6]{CP}.
Given a rank-two permutation matrix $u_\sigma$ we define a selfadjoint subspace $\Xi$ of
$\F_2^1 := {\rm span}\{s_i s_j^*, 1 \leq i,j \leq 2\} \simeq M_2$
as $\bigcap_{s \in {\mathbb N}} (K^*)^s \F_2^1 K^s$,
where $K = {\rm span} \{u_\sigma s_1, u_\sigma s_2 \} \subset \O_2$.
We will need the next result, that gathers together few facts from \cite[Section 6]{CP}.

\begin{theorem}\label{CP}
Let $u$ be a unitary in $ \F_2^2 := M_2 \otimes M_2$ and assume that the corresponding subspace
$\Xi$ satisfies $\Xi^2 \subseteq \Xi$.
Suppose in addition that 
\begin{itemize}
\item[(a)] $\omega(a\lambda_u(b)) = \omega(a)\omega(b)$ for all $a \in \Xi$, $b \in \F_2^1$
\item[(b)] $E_1 (u^* a u) = \omega(a) 1$ for every $a \in \Xi$, where $E_1: \pi_\omega(\O_2)'' \to \F_2^1$
is the $\omega$-invariant conditional expectation. 
\end{itemize}
Then
$${\rm Ind}(\lambda_u) = \dim(\Xi) \ . $$ 
\end{theorem}

\begin{proof}
It follows directly from \cite[Proposition 6.1]{CP} and \cite[Theorem 6.5 (c)]{CP} (or \cite[Corollary 6.6]{CP}).
\end{proof}

\medskip
\noindent{\bf The transformation induced by $\sigma_{12}$.}
One can check that the subspace $\Xi \subset M_2$ is the linear span of
$1$ and $s_1 s_2^* + s_2 s_1^*$.
Moreover, the conditions (a) and (b) in Theorem \ref{CP}
are easily verified and 
thus
${\rm Ind}(\rho_{12}) = 2$.

\medskip
\noindent{\bf The transformation induced by $\sigma_{13}$.}
One can check that the subspace $\Xi \subset M_2$ is the linear span of
$s_1 s_1^*$ and $s_2 s_2^*$.
Moreover, the conditions (a) and (b) in Theorem \ref{CP}
are easily verified and thus
${\rm Ind}(\rho_{13}) = 2$.

\medskip
\noindent{\bf The transformation induced by $\sigma_{142}$.}
One can check that in this case $\Xi = \F_1^1$. 
Moreover, the conditions (a) and (b) in Theorem \cite{CP}
are easily verified and thus ${\rm Ind}(\rho_{142}) = 4$.
In alternative,
one can apply \cite[Corollary 7.6]{CP} with $U = u_{142}$ and $V = -u_{134}$.
We omit the lengthy but straightforward calculation.
Notice that the restriction of $\rho_{142}$ to $\F_2$ can easily be shown to be reducible.

\medskip
\noindent{\bf The remaining transformations.}
One has $\rho_{34} = \rho_{(12)(34)} \circ \rho_{12}$,
$\rho_{1324} = {\rm Ad}(u_{(13)(24)}) \circ \rho_{12}$,
$\rho_{1423} = {\rm Ad}(u_{(13)(24)}) \circ \rho_{34}$,
$\rho_{24} = \rho_{13} \circ \rho_{(13)(24)}$,
$\rho_{1234} = {\rm Ad}(u_{(13)(24)}) \circ \rho_{24}$
and $\rho_{1432} = {\rm Ad}(u_{(13)(24)}) \circ \rho_{13}$
and thus, by multiplicativity,
${\rm Ind}(\rho_{34}) = {\rm Ind}(\rho_{1324}) = {\rm Ind}(\rho_{1423}) =
{\rm Ind}(\rho_{24}) = {\rm Ind}(\rho_{1234}) = {\rm Ind}(\rho_{1432}) = 2$.
Also,
one has $\rho_{134} = \rho_{142} \circ \rho_{(13)(24)}$,
thus ${\rm Ind}(\rho_{134})={\rm Ind}(\rho_{142})=4$.

{\footnotesize
\[
\begin{array}{lcccccc}
\multicolumn{6}{c}{\mbox{Table 1. Entropy and index of the `rank 2' permutation endomorphisms of $\O_2$.}}\\
\\
\hline  \rho_{\sigma}& \rho_{\sigma}(s_{1})& \rho_{\sigma}(s_{2})& property &
{\rm ht}(\rho_\sigma)& {\rm ht}(\rho_\sigma|_{\D_2}) & {\rm Ind}(\rho_\sigma)\\
\hline
\rho_{id} = {\rm id} & s_{1} & s_{2}& {\rm inn} & 0&0 & 1\\
\rho_{12}   &s_{12,1}+s_{11,2}& s_{2}& {\rm irr} & \log 2& 0 & 2\\
\rho_{13}   &s_{21,1}+s_{12,2}& s_{11,1}+s_{22,2} & {\rm irr}& \log 2&\log 2 & 2\\
\rho_{14}   &s_{22,1}+s_{12,2}& s_{21,1}+s_{11,2}& {\rm red}& \log 2&\log 2 & 4\\
\rho_{23}   &s_{11,1}+s_{21,2}& s_{12,1}+s_{22,2}& {\rm red}& \log 2&\log 2 & 4\\
\rho_{24}   &s_{11,1}+s_{22,2}&s_{21,1}+s_{12,2}& {\rm irr}& \log 2&\log 2 & 2\\
\rho_{34}   &s_{1}& s_{22,1}+s_{21,2}& {\rm irr}&\log 2&0 & 2\\
\rho_{123} &s_{12,1}+s_{21,2}& s_{11,1}+s_{22,2}& {\rm red}&\log 2&\log 2 & 4\\
\rho_{132} &s_{21,1}+s_{11,2}& s_{12,1}+s_{22,2}& {\rm red}&\log 2&\log 2 & 4\\
\rho_{124} &s_{12,1}+s_{22,2}& s_{21,1}+s_{11,2}& {\rm red}&\log 2&\log 2 & 4\\
\rho_{142} &s_{22,1}+s_{11,2}& s_{21,1}+s_{12,2}& {\rm irr}&\log 2&\log 2 & 4\\
\rho_{134} (\simeq \rho_{142}) &s_{21,1}+s_{12,2}& s_{22,1}+s_{11,2}& {\rm irr}&\log 2&\log 2 & 4\\
\rho_{143} &s_{22,1}+s_{12,2}& s_{11,1}+s_{21,2}& {\rm red}&\log 2&\log 2 & 4\\
\rho_{234} &s_{11,1}+s_{21,2}& s_{22,1}+s_{12,2}& {\rm red}&\log 2&\log 2 & 4\\
\rho_{243} (\simeq \rho_{123}) &s_{11,1}+s_{22,2}& s_{12,1}+s_{21,2}& {\rm red}&\log 2&\log 2 & 4\\
\rho_{1234} (\simeq \rho_{24}) &s_{12,1}+s_{21,2}& s_{22,1}+s_{11,2}& {\rm irr}&\log 2&\log 2 & 2\\
\rho_{1243} &s_{12,1}+s_{22,2}& s_{11,1}+s_{21,2}& {\rm red}&\log 2&\log 2 & 4\\
\rho_{1324} (\simeq \rho_{12}) &s_{2}& s_{12,1}+s_{11,2}& {\rm irr}&\log 2&0 & 2\\
\rho_{1342} &s_{21,1}+s_{11,2}& s_{22,1}+s_{12,2}& {\rm red}&\log 2&\log 2 & 4\\
\rho_{1423} (\simeq \rho_{34}) &s_{22,1}+s_{21,2}& s_{1}& {\rm irr}&\log 2&0 & 2\\
\rho_{1432} (\simeq \rho_{13}) &s_{22,1}+s_{11,2}& s_{12,1}+s_{21,2}& {\rm irr}&\log 2&\log 2 & 2\\
\rho_{(12)(34)} (\simeq \rho_{(13)(24)}) &s_{12,1}+s_{11,2}& s_{22,1}+s_{21,2}& {\rm out}& 0&0 & 1\\
\rho_{(13)(24)} &s_{2}& s_{1}& {\rm out}&0&0 & 1\\
\rho_{(14)(23)} (\simeq {\rm id}) &s_{22,1}+s_{21,2}& s_{12,1}+s_{11,2}& {\rm inn}&0&0 & 1\\
\hline
\end{array}
\]
}

\section{Further comments}
In Table 1 above we collected the currently known key characteristics
of the rank-two permutation endomorphisms of ${\mathcal O}_2$.\footnote{In the fourth
column, inn, irr, red and out stand for ``inner automorphism'', ``irreducible endomorphism'',
``reducible endomorphism'' and ``outer automorphism'', respectively, as endomorphisms of $\O_2$.}
At this point we do not know of any general direct relation between their indices (the
seventh column of the table) and entropies (the fifth and sixth columns).
However, it is natural
to wonder whether there is a more direct connection between these index computations
and the framework in \cite{SZ}. Certainly the (minimal) index is related to the Pimsner-Popa entropy \cite{PP}.
From this point of view, it is quite intriguing that
the inequality in \cite[Theorem 2.1]{SZ} is completely analogous to the one already known for the index \cite{l}
(see also \cite{CP}, bottom of p.373),
thus suggesting the existence of some relation between
(square root of)
index and topological entropy, cf. \cite[Section 10]{St}.
Notice that the index pertains to von Neumann algebras while the topological entropy refers to $C^*$-algebras,
so it would probably make sense also to investigate Watatani's notion of index \cite{Wa} in a purely $C^*$-framework
to spell out more details.

It is also worth pointing out that in the situation under consideration
both the quantities ${\rm ht}(\rho_\sigma)$ and ${\rm ht}(\rho_\sigma|_{\D_2})$ studied in \cite{SZ}
attain only two values, 0 and $\log 2$.
Which value occurs depends on whether the given endomorphism or its restriction
to the diagonal is actually an automorphism or not,
a fact that certainly deserves more thorough investigations.
Indeed, as shown in \cite{SZ}, ${\rm ht}(\rho_\sigma)=0$ if and only if $\rho_\sigma$ is an automorphism of $\O_2$.
Furthermore, the four proper endomorphisms $\rho_\sigma$ with ${\rm ht}(\rho_\sigma|_{\D_2})=0$ are precisely those
four proper endomorphisms which restrict to automorphisms of the diagonal $\D_2$.\footnote{The irreducibility
of such endomorphisms is automatic \cite[Proposition 1.1]{CKS}.}
This is easily verified by the method developed in \cite{CS,Sz} (see Theorem 3.4, Lemma 4.5, Corollary 4.6 and Subsection 5.1 in \cite{CS}).
Also note that the four automorphisms of $\O_2$ appearing in Table 1 all are of finite order
(order two) and thus their topological entropy is 0 \cite[Proposition 4.2]{Vo}.
It is shown in \cite[Subsection 5.2]{CS} that all rank-three automorphisms of $\O_2$
(and there are exactly 48 of them) have also finite orders.
The first infinite order automorphisms appear in rank-four \cite[Subsection 5.3]{CS}
and it would be interesting to determine their topological entropies.
In the case of $\O_n$ with $n \geq 3$ infinite order automorphisms appear already at rank two
and the question of determining their entropies remains open too.

\bigskip

Regarding index computations, for a fixed value of $n$, say $n=2$ in our case, the next step of complexity would be to
look at the $2^3 ! $ rank-three permutation endomorphisms.
Among those, the ones providing genuine automorphisms of $\O_2$ or $\D_2$ have been completely classified in \cite{CS}.
Still, one can not exclude the possibility that some proper endomorphisms of $\O_2$ extend
to automorphisms of the associated factor.
Some examples of this kind were announced for the case $n=3$ and $k=2$ in \cite{Ak}.


\smallskip\noindent
Roberto Conti\\
Mathematics, School of Mathematical and Physical Sciences \\
The University of Newcastle, Callaghan, NSW 2308, Australia
\\ E-mail: Roberto.Conti@newcastle.edu.au \\

\smallskip \noindent
Wojciech Szyma{\'n}ski\\
Mathematics, School of Mathematical and Physical Sciences \\
The University of Newcastle, Callaghan, NSW 2308, Australia
\\ E-mail: Wojciech.Szymanski@newcastle.edu.au \\

\noindent Present address: \\
Department of Mathematics and Computer Science \\
The University of Southern Denmark \\
Campusvej 55, DK-5230 Odense M, Denmark \\
E-mail: szymanski@imada.sdu.dk


\begin{thebibliography}{99}

\bibitem{Ak} P. T. Akemann,
{\em On a class of endomorphisms of the hyperfinite $II_1$ factor},
Doctoral Dissertation, UC Berkeley, 1997.

\bibitem{Ch} M. Choda,
{\em Entropy of Cuntz's canonical endomorphism},
Pacific J. Math. {\bf 190} (1999), 235--245.

\bibitem{CF} R. Conti and F. Fidaleo, {\it Braided endomorphisms of Cuntz algebras},
Math. Scand. {\bf 87} (2000), 93--114.

\bibitem{CP} R. Conti and C. Pinzari, {\it Remarks on endomorphisms of Cuntz
algebras}, J. Funct. Anal. {\bf 142} (1996), 369--405.

\bibitem{CS} R. Conti and W. Szyma{\'n}ski,
{\em Labeled trees and localized automorphisms of the Cuntz algebras}, submitted,
arXiv:0805.4654.

\bibitem{CKS} R. Conti, J. Kimberley and W. Szyma{\'n}ski,
{\em More localized automorphisms of the Cuntz algebras}, submitted,
arXiv:0808.2843.

\bibitem{Cun1} J. Cuntz, {\it Simple $C^*$-algebras generated by isometries},
Commun. Math. Phys. {\bf 57} (1977), 173--185.

\bibitem{Cun2} J. Cuntz, {\it Automorphisms of certain simple $C^*$-algebras},
in {\it Quantum fields-algebras-processes}, ed. L. Streit, Springer 1980.

\bibitem{Go} R. Gohm, {\it A probabilistic index for completely positive maps
and an application}, J. Operator Theory {\bf 54} (2005), 339--361.

\bibitem{Hi} F. Hiai,
{\em Minimizing indices of conditional expectations onto a subfactor},
Publ. Res. Inst. Math. Sci., Kyoto Univ. {\bf 24} (1988), 673--678.

\bibitem{i} M. Izumi, {\it Subalgebras of infinite $C^*$-algebras with
finite Watatani indices. I. Cuntz algebras}, Commun. Math. Phys. {\bf 155}
(1993), 157--182.

\bibitem{Jo} V. F. R. Jones,
{\em Index for subfactors},
Invent. Math. {\bf 72} (1983), 1--25.

\bibitem{j0} V. F. R. Jones,
{\em Subfactors and knots}. CBMS Regional Conference Series in Mathematics, 80.
Published for the Conference Board of the Mathematical Sciences, Washington, DC;
by the American Mathematical Society, Providence, RI, 1991.

\bibitem{j} V. F. R. Jones,
{\em On a family of almost commuting endomorphisms},
J. Funct. Anal. {\bf 122} (1994), 84--90.

\bibitem{JS} V. F. R. Jones and V. S. Sunder, {\em Introduction to subfactors},
London Math. Soc. Lecture Note Ser. {\bf 234}, Cambridge University Press,
Cambridge, 1997.

\bibitem{KaLo} Y. Kawahigashi and R. Longo,
{\em Classification of Local Conformal Nets. Case $c < 1$},
Ann. of Math. {\bf 160} (2004), 493--522.

\bibitem{Kaw0} K. Kawamura,
{\em Algebra of sectors}, RIMS-1450 (2004).

\bibitem{Kaw1} K. Kawamura,
{\em Polynomial endomorphisms of the Cuntz algebras arising from permutations. I.
General theory},  Lett. Math. Phys.  {\bf 71}  (2005),  149--158.

\bibitem{Kaw2} K. Kawamura,
{\em Branching laws for polynomial endomorphisms of Cuntz algebras arising from
permutations},  Lett. Math. Phys.  {\bf 77}  (2006),   111--126.

\bibitem{Ko} H. Kosaki,
{\em Extension of Jones' theory on index to arbitrary factors},
J. Funct. Anal. {\bf 66} (1986), 123--140.

\bibitem{Lo1} R. Longo,
{\em Index of subfactors and statistics of quantum fields. I},
Commun. Math. Phys. {\bf 126} (1989), 217--247.

\bibitem{Lo2} R. Longo,
{\em Index of subfactors and statistics of quantum fields. II.
Correspondences, braid group statistics and Jones polynomial},
Commun. Math. Phys. {\bf 130} (1990), 285--309.

\bibitem{l} R. Longo, {\em A duality for Hopf algebras and for subfactors. I},
Commun. Math. Phys. {\bf 159} (1994), 133--150.

\bibitem{l2} R. Longo, {\em Notes for a quantum index theorem},
Commun. Math. Phys. {\bf 222} (2001), 45--96.

\bibitem{PP} M. Pimsner and S. Popa, {\em Entropy and index for subfactors},
Ann. Sci. \'Ecole Norm. Sup. {\bf 19} (1986), 57--106.

\bibitem{SZ} A. Skalski and J. Zacharias,
{\em Noncommutative topological entropy of endomorphisms of Cuntz algebras},
to appear in {\em Lett. Math. Phys.}, arXiv:0804.4373.

\bibitem{St} E. St\o rmer,
{\em A survey of noncommutative dynamical entropy},
Classification of nuclear $C\sp *$-algebras. Entropy in operator algebras, 147--198,
Encyclopaedia Math. Sci., 126, Springer, Berlin, 2002.

\bibitem{Sz} W. Szyma{\'n}ski,
{\it On localized automorphisms of the Cuntz algebras which preserve the diagonal subalgebra},
in `New Development of Operator Algebras', R.I.M.S. K\^{o}ky\^{u}roku {\bf 1587} (2008), 109--115.

\bibitem{Vo} D. V. Voiculescu,
{\em Dynamical approximation entropies and topological entropy in operator algebras},
Commun. Math. Phys. {\bf 170} (1995), 249--281.

\bibitem{Wa} Y. Watatani, {\em Index for $C^*$-subalgebras},
Mem. Amer. Math. Soc. 83 (1990), no. 424.

\end{thebibliography}
\end{document}